\documentclass{article}
\usepackage[utf8]{inputenc}
\usepackage{listings}
\usepackage{amsmath}
\usepackage{amsthm}
\usepackage{amsfonts}
\usepackage{amssymb}
\usepackage[english]{babel}
\usepackage[a4paper]{geometry}
\usepackage{graphicx}
\usepackage{color,psfrag}
\usepackage{fancyhdr}
\usepackage{cite}
\usepackage{makeidx}
\usepackage{enumerate}
\usepackage{mathpazo}
\usepackage{algorithm}
\usepackage{algorithmic}
\usepackage{mathrsfs}
\usepackage{float}
\usepackage{datetime}
\usepackage{titlesec}
\usepackage[labelfont=bf,textfont=it,labelsep=period,width=.8\textwidth]{caption}
\usepackage[unicode,pdfstartview={XYZ null null 1},pdfview={XYZ null null 1},pdfpagemode=UseOutlines]{hyperref}
\DeclareFontFamily{U}{wncy}{}
\DeclareFontShape{U}{wncy}{m}{n}{<->wncyr10}{}
\DeclareSymbolFont{mcy}{U}{wncy}{m}{n}
\DeclareMathSymbol{\Sh}{\mathord}{mcy}{"58}

\newcommand{\A}{{\cal A}}

\def\B{{\cal B}}

\newcommand{\R}{\mathbb R}

\newdate{date}{10}{03}{2017}

\newcommand{\F}{\mathbb F}

\newtheorem{theorem}{Theorem}

\newtheorem{lemma}[theorem]{Lemma}

\newtheorem{proposition}[theorem]{Proposition}
\newtheorem{protocol}{Protocol}
\newtheorem{definition}[theorem]{Definition}

\newtheorem{heuristic}{Heuristic}
\newtheorem{problem}{Problem}
\theoremstyle{remark}
\theoremstyle{remark}\newtheorem{remark}[theorem]{Remark}

\title{An identification system based on the explicit isomorphism problem}
\author{\normalsize
 \begin{minipage}{0.3\linewidth}
    \large
    Sándor Z. Kiss \\
    \footnotesize
    Department of Algebra, Budapest
    University of Technology and Economics \\
    \texttt{kisspest@cs.elte.hu} \\
    \normalsize
  \end{minipage}
  \qquad
 \begin{minipage}{0.3\linewidth}
    \large
    P\'eter Kutas \\
    \footnotesize
School of Computer Science, University of Birmingham\\
    \texttt{p.kutas@bham.ac.uk} \\
    \normalsize
  \end{minipage}
}

\begin{document}

\maketitle
\begin{abstract}
\noindent We propose a new identification system based on algorithmic problems related to computing isomorphisms between central simple algebras. We design a statistical zero knowledge protocol which relies on the hardness of computing isomorphisms between orders in division algebras which generalizes a protocol by Hartung and Schnorr, which relies on the hardness of integral equivalence of quadratic forms.
\end{abstract}


\bigskip
\noindent

\section{Introduction}

In this paper we propose an identification system based on an algorithmic problem related to the following problem from computational algebra. Let $\A$ be a finite-dimensional associative algebra over a field $\mathbb{K}$.  Let $b_1,\dots,b_{m}$ be a basis of $\A$. Then the products $b_ib_j$ can be expressed as linear combinations of the basis elements: $b_ib_j=\sum_{k=1}^{m} \gamma_{ijk}b_k.$
The $\gamma_{ijk}$s are called structure constants and we consider $\A$ to be given by a collection of structure constants. Assume that $\A$ is isomorphic to $M_n(\mathbb{K})$, the algebra of $n$ times $n$ matrices over $\mathbb{K}$. The algorithmic task is to compute an isomorphism between $\A$ and $M_n(\mathbb{K})$. We will refer to this problem as the explicit isomorphism problem. 

This is a well studied problem in computational algebra \cite{Cremona3},\cite{IKR},\cite{IRS}, \cite{K}, \cite{Pilnikova}. It has connections to arithmetic geometry \cite{Cremona1},\cite{Cremona2}, \cite{Cremona3}, norm equations \cite{IRS}, parametrization of algebraic varieties \cite{de Graaf} and error-correcting codes \cite{GTLN}. The best known algorithm in the $\mathbb{K}=\mathbb{Q}$ case is due to Ivanyos, Rónyai and Schicho \cite{IRS}. The algorithm uses an oracle for integer factorization and the running time of the algorithm is polynomial in the size of the structure constants but is exponential in the $n$ (the degree of the matrix algebra) which implies that it is only practical for very small $n$. The algorithm of \cite{IRS} can also be used to compute isomorphisms between division algebras by a reduction to the original explicit isomorphism problem. This reduction on the other hand comes at the cost of squaring the dimension. To our knowledge the difficulty of this problem has never been exploited for cryptographic purposes. 

Hartung and Schnorr \cite{Hartung} proposed an identification system which relies on the difficulty of finding an explicit equivalence of integral quadratic forms.
In a sense our scheme can be thought of as a higher degree generalization of the protocol in \cite{Hartung} as the equivalence problem of rational quadratic forms is similar to the isomorphism problem of quaternion algebras. 

We introduce new computational problems (e.g., the order isomorphism problem) which are naturally harder problems than the above discussed explicit isomorphism problem. The complexity of them is unclear at the moment but there is some evidence that due to their "integral" nature they are indeed harder.

The paper is organized as follows. In Section 2 we summarize all known results and computational assumptions which we will use later on. In Section 3 we give the detailed description of our protocol and provide security proofs. Finally, in Section 4 we give a toy example of the protocol described in Section 3.

\section{Preliminaries}
\subsection{Theoretical background}
In this subsection we give a brief overview of the theoretical results needed for the description of our protocol. The reader is referred to \cite[Chapter 12]{Pi} on facts about central simple algebras. 

\begin{definition}[\cite{Pi}, p.44]
A nonzero associative algebra $\A$ over a field $\mathbb{K}$ is simple if it has no nontrivial two-sided ideals. 
\end{definition}
\noindent It is a well-known theorem of Wedderburn that every finite dimensional simple algebra over a field $\mathbb{K}$ is isomorphic to a full matrix algebra over some division algebra whose center is an extension of $\mathbb{K}$.
\begin{definition}[\cite{Pi}, p.224]
A simple algebra $\A$ is called central simple over $\mathbb{K}$ if its center is exactly $\mathbb{K}$.
\end{definition}
\noindent The tensor product of two finite-dimensional central simple $\mathbb{K}$-algebras is again a central simple $\mathbb{K}$-algebra. Two central simple algebras over $\mathbb{K}$ are Brauer equivalent if their underlying division algebras are isomorphic (note that by Wedderburn's theorem, a  central simple algebra is a full matrix algebra over a division algebra). Equivalence classes of central simple algebras form a group under the tensor product, called the Brauer group of the field $\mathbb{K}$. This implies that in order to understand central simple algebras over a fixed field one has to understand the division algebras over that field.

\begin{definition}[\cite{Pi}, p.277]
Let $\mathbb{K}$ be a field and let $\mathbb{L}$ be a cyclic extension of $\mathbb{K}$ (i.e., a Galois extension whose Galois group is cyclic) of degree $n$ . Let $\sigma$ be a generator of the Galois group. Let $a\in \mathbb{K}$. Then the following algebra $\A$ is called a cyclic algebra:
\begin{enumerate}
    \item $u^n=a\cdot 1$
    \item $\A=\oplus_{i=0}^{n-1}\mathbb{L}u^i$
    \item $u^{-1}lu=\sigma(l)$ for every $l\in \mathbb{L}$
\end{enumerate}
This algebra is denoted by $(\mathbb{L}|\mathbb{K},\sigma, a)$. 
\end{definition}
\noindent It is well known (see \cite[Chapter 15]{Pi}) that a cyclic algebra is a central simple algebra over $\mathbb{K}$ of dimension $n^2$. Moreover, the  following is true: 

\begin{theorem}\label{cyclic}(\cite{Pi}, p.278)
Let $\mathbb{L}$ be a cyclic extension of $\mathbb{K}$ and let $a\in \mathbb{K}\setminus\{0\}$. A cyclic algebra $(\mathbb{L}|\mathbb{K},\sigma, a)$ is isomorphic to $M_n(\mathbb{K})$ if and only if $a$ is a norm in the extension $\mathbb{L}|\mathbb{K}$. 
\end{theorem}
\noindent We define orders in central simple algebras:
\begin{definition}[\cite{Reiner}, p.108]
Let $R$ be an integral domain with quotient field $\mathbb{K}$. Let $\A$ be a central simple $\mathbb{K}$-algebra. A subring $O$ of $\A$ is an order if it contains 1 and is a finitely generated $R$-module which contains a $\mathbb{K}$-basis of $\A$ (i.e., $O\otimes_R \mathbb{K}=\A$).
\end{definition}
\noindent An order is called maximal, if it is maximal with respect to inclusion. Maximal orders are non-commutative analogues of the ring of integers in algebraic number fields. For further details on maximal orders the reader is referred to Reiner's monograph \cite{Reiner}.
\begin{theorem}\label{NS}(Noether-Skolem \cite{Pi}, p.230)
Let $\A$ be a finite dimensional central simple $\mathbb{K}$-algebra and let $\B$ be a simple $\mathbb{K}$-algebra. Let $f,g: \B \rightarrow \A$ be two $\mathbb{K}$-algebra homomorphisms. Then there exists an invertible element $x \in \A$, such that $f(b) = xg(b)x^{-1}$ for all $b \in \B$.
\end{theorem}

\subsection{Algorithmic background and computational assumptions}

\subsubsection{Known results}
In this subsection we give a brief overview of the algorithmic history of the explicit isomorphism problem.

Let $\A$ be an associative algebra given by a collection of structure constants. It is a natural algorithmic problem to compute the structure of $\A$, i.e., compute its Jacobson radical $rad~ \A$, compute the Wedderburn decomposition of $\A/rad~ \A$ and finally compute an explicit isomorphism between the simple components of $\A/rad~ \A$ and $M_{n_i}(\mathcal{D}_i)$ where the $\mathcal{D}_i$s are division algebras over $\mathbb{K}$ and $M_{n_i}(\mathcal{D}_i)$ denotes the algebra of $n_i\times n_i$ matrices over $\mathcal{D}_i$. The problem has been studied for various fields $\mathbb{K}$, including finite fields, the field of complex and real numbers, global function fields and algebraic number fields. There exists a polynomial-time algorithm for computing the radical of $\A$ over these fields \cite{CIW}. There also exist efficient algorithms for every task over finite fields \cite{FR},\cite{Ro1} and the field of real and complex numbers \cite{Eb3}. Finally, when $\mathbb{K}=\mathbb{F}_q(t)$, the field of rational functions over a finite field $\F_q$, then there exist efficient algorithms for computing Wedderburn decompositions \cite{IRSz} and for computing explicit isomorphisms between full matrix algebras over $\mathbb{F}_q(t)$ \cite{IKR}. 

The case when $\mathbb{K}=\mathbb{Q}$ is particularly interesting due to its applicability to various algorithmic problems. Wedderburn decomposition can again be achieved in polynomial time \cite{FR}, but computing isomorphisms between central simple $\mathbb{K}$-algebras is much harder. Rónyai \cite{Ronyai} showed that computing an explicit isomorphism between $\A$ (given by structure constants) and $M_2(\mathbb{Q})$ is at least as hard as factoring integers. On the other hand, Ivanyos, Rónyai and Schicho \cite{IRS} proposed an algorithm to compute an isomorphism between $\A$ and $M_n(\mathbb{Q})$ which is allowed to call an oracle for factoring integers. The running time of the algorithm is polynomial in the size of the structure constants, but it is exponential in $n$. More precisely, let $m=n^2$ and let $$c_m=\gamma_m^{\frac{m}{2}}\left( \frac{3}{2}\right )^m2^{\frac{m(m-1)}{2}},$$
where $\gamma_{m}$ is Hermite's constant.
The last step of the algorithm from \cite{IRS} generates roughly $c_m^m$ linear combinations and checks whether any of them has rank 1 as a matrix. This number is independent of the size of the structure constants, so when $m$ is bounded, the algorithm is technically a polynomial-time algorithm. However, when $n\geq 5$ the search space is too large for this computation to be achievable in reasonable amount of time. Making this algorithm practical would be of immense number theoretical intereset as it would speed up $n$-descent of elliptic curves \cite{Cremona3} which is one of the most promising techniques for computing generators of the Mordell-Weil group of elliptic curves.

In \cite{IRS} the authors also study the isomorphism problem of division algebras. They reduce the problem of finding explicit isomorphisms between division algebras of degree $n$ over a number field $\mathbb{K}$ to the explicit isomorphism problem between an algebra $\A$ and $M_{n^2}(\mathbb{K})$ \cite[Section 4]{IRS}. Note that this suggests that the isomorphism problem of division algebras is harder as it requires the solution of an explicit isomorphism problem for full matrix algebras of degree $n^2$ as opposed to $n$. 

If $\A$ is given by a cyclic algebra presentation, then finding an isomorphism between $\A$ and $M_n(\mathbb{Q})$ is equivalent to solving a norm equation over a cyclic extension which is a classical hard problem in computational number theory. Furthermore, there is no known polynomial-time algorithm for computing a cyclic algebra presentation from a structure constant representation when $n\geq 5$. So it seems that the explicit isomorphism problem is harder than solving norm equations in cyclic extensions. 

It is a natural question to study related isomorphisms problems, such as the isomorphism problem of orders in division algebras. This has not extensively been studied but there is some evidence that this problem is harder than the previous problems considered. First, note that if one can construct an isomorphism between two orders, then that extends to an isomorphism of the underlying algebras (computing an order in an algebra is easy, one just multiplies the basis elements with a suitable integer to make structure constants integral). The relation of order isomorphism and algebra isomorphism is similar to the relation between finding integer solutions and rational solutions to a diophantine equation. The problem of equivalence of quadratic forms is studied in \cite{Hartung} where they show that rational equivalence can be computed by a polynomial-time algorithm which is allowed to call an oracle for factoring integers. On the other hand, they also show that the problem of integral equivalence is NP-hard. This is similar as the relation between isomorphisms of division algebras and orders. This provides some evidence that the order isomorphism problem is harder than, and extends the isomorphism problem of algebras. 

\subsubsection{Computational assumptions}

We list hard problems and list our computational assumptions which are needed for our scheme. We start by restating the explicit isomorphism problem:

\begin{problem}\label{EIP}
Let $\A$ be an algebra isomorphic to $M_n(\mathbb{Q})$ given by structure constants. The \emph{Explicit Isomorphism Problem (EIP)} is to find an isomorphism between $\A$ and $M_n(\mathbb{Q})$. 
\end{problem}
\noindent In order to be able to consider more general problems, we formalize isomorphism problems in such a way that checking if a map is really and algebra isomorphism can be accomplished efficiently. First we give a slightly different interpretation of the structure constant representation.
Let $\A$ be an algebra of dimension $m$ over a field $\mathbb{K}$. Then multiplication from the left by any element $a\in\A$ is a $\mathbb{K}$-linear map, thus can be described by an $m\times m$ matrix with entries from $\mathbb{K}$. This provides an embedding of $\A$ into the full matrix algebra $M_m(\mathbb{K})$. This representation of the algebra is called the \emph{(left-)regular representation}. It is clear that specifying an algebra by a collection of structure constants is exactly the same as providing its regular representation. 
An isomorphism between algebras given by structure constants can be specified in various ways. One way is to describe it as a vector space isomorphism. However, checking that it is also ring isomorhism might be costly. Instead we consider both algebras given by their regular representation and then an isomorphism can be described as conjugation by a suitable matrix.

\begin{problem}\label{DAIP}
Let $\A,\B$ be isomorphic division algebras of dimension $n^2$ over $\mathbb{Q}$ given by their regular representation. The \emph{Division Algebra Isomorphism Problem} is to find an invertible matrix $M\in M_{n^2}(\mathbb{Q})$ such that $\B=M^{-1}\A M$. 
\end{problem}
\noindent The Noether-Skolem theorem implies the existence of a suitable $M$. Now we state the order isomorphism problem: 
\begin{problem}\label{OIP}
Let $\A,\B$ be isomorphic division algebras of dimension $n^2$ over $\mathbb{Q}$ and let $\Gamma_{\A},\Gamma_{\B}$ be orders in $\A$ and $\B$ respectively. The \emph{Order Isomorphism Problem} is to find an invertible matrix $M\in M_{n^2}(\mathbb{Q})$ such that $\Gamma_{\B}=M^{-1}\Gamma_{\A}M$. 
\end{problem}
\noindent Finally we consider a slightly more general version of the order isomorphism problem when we require the conjugating matrix to be integral. 
\begin{problem}\label{IOIP}
Let $\A,\B$ be isomorphic division algebras of dimension $n^2$ over $\mathbb{Q}$ and let $\Gamma_{\A},\Gamma_{\B}$ be orders in $\A$ and $\B$ respectively. The \emph{Integer Order Isomorphism Problem} is to find an invertible matrix $M\in M_{n^2}(\mathbb{Z})$ with $\det(M) = \pm 1$ such that $\Gamma_{\B}=M^{-1}\Gamma_{\A}M$. 
\end{problem}

\noindent Analogies suggest that Problem \ref{IOIP} is the hardest amongst these problems. An example of Problem \ref{IOIP} for central simple algebras of degree two (i.e., quaternion algebras) is given in Section \ref{sec:example} as part of the toy example of our protocol. 

In \cite{IRS} it is shown that if two central simple algebras are isomorphic then there exists an isomorphism that can be represented by a matrix of polynomial size. However, it is not obvious that the statement also holds for orders. This motivates the definition of the following algorithmic problem: 
\begin{problem}\label{IOIP2}
Let $\A,\B$ be isomorphic division algebras of dimension $n^2$ over $\mathbb{Q}$ and let $\Gamma_{\A},\Gamma_{\B}$ be orders in $\A$ and $\B$ respectively. The \emph{Integer Order Isomorphism Problem With Restricted Coefficients} is to find an invertible matrix $M\in M_{n^2}(\mathbb{Z})$ with $\det(M) = \pm 1$ such that $\Gamma_{\B}=M^{-1}\Gamma_{\A}M$ and every entry $a_{i,j}$ of $M$ one has that $|a_{i,j}|<t$ for some constant $t$. 
\end{problem}
\noindent Our main protocol will rely on the hardness of this problem.
\begin{lemma}\label{lem:np}
Problem \ref{IOIP2} is in NP.
\end{lemma}
\begin{proof}
Let $\Gamma_{\A}$ and $\Gamma_{\B}$ be orders such that $\Gamma_{\B}=M^{-1}\Gamma_{\A} M$ and $M\in M_{n^2}(\mathbb{Z})$ with $\det(M) = \pm 1$. Then we show that $M$ is a polynomial-time witness for Problem \ref{IOIP2}. First, one checks that $M$ has integer coefficients which are bounded in absolute value by $t$ and has determinant equal to $\pm 1$. Then one computes the $M^{-1}bM$ for every basis element and computes the coefficient of the above matrices in this basis by solving a system of linear equations. If the coefficients are integers and the transition matrix has determinant $\pm 1$, then $M$ indeed induces an isomorphism. In other words, one has to check that the $\mathbb{Z}$-lattice generated by the elements $M^{-1}bM$ is equal to $\Gamma_{\A}$ which is an instance of equality of lattices. 
\end{proof}

\subsection{Interactive Proof systems}
In this section we give a short survey about interactive proof systems and zero knowledge protocols. We follow \cite{DelfsKnebl}, which contains an excellent summary about identification schemes.
An interactive proof system consists of two participants: the prover and the verifier. The aim of the prover is to convince the verifier that he knows some secret information, which is called prover's secret. During the whole process the prover and the verifier send and receive messages from each other and both of them perform some computations. The communication between the prover and the verifier consists of challenges by the verifier and responses by the prover. In general, the prover begins and the verifier finishes the protocol. The verifier accepts or rejects depending on the prover's answers to all of the verifier's challenges.  

For an interactive proof system, there are the following two requirements. (\cite{DelfsKnebl}, p.118.)
\begin{enumerate}
\item[1.] \textbf{Completeness.} If the prover knows the prover's secret, then the verifier will always accept the prover's proof.
\item[2.] \textbf{Soundness.} If the prover can convince the verifier with reasonable probability, then he knows the prover's secret.
\end{enumerate}
A prover or a verifier is called honest prover or honest verifier if he follows the steps specified in the protocol; otherwise he is called dishonest or fraudulent prover or verifier. 

Now, we give a formal definition for the zero-knowledge property (see \cite{DelfsKnebl}, subsection 4.2.3.) We denote the algorithm of the honest prover by $P$, the algorithm of the honest verifier by $V$, and the algorithm of an arbitrary verifier by $V^{*}$. Note that $V^{*}$ can be a fraudulent verifier. We denote an interactive proof system, including the interaction between $P$ and $V$ by $(P,V)$. 
Assume that the interactive proof takes $n$ steps. 
In each step a message is sent, and we can assume that the prover starts with the first step. Let $m_{1}, m_{3}, \dots{}$ be the messages sent from the prover to the verifier and let $m_{2}, m_{4}, \dots{}$ be the messages sent from the verifier to the prover, where $m_{i}$ denotes the message sent in the $i$-th step.
We define the transcript of the joint computation of $P$ and $V^{*}$ 
the common input $x$ of $(P,V)$ by $tr_{P,V^{*}}(x) = (m_{1}, m_{2}, \dots{}, m_{n})$, where $tr_{P,V^{*}}(x)$ is called an accepting transcript if $V^{*}$ accepts after the last step. 

For a given verifier $V^{*}$, an algorithm $S$ which generates valid accepting transcripts for $(P,V^{*})$ without communicating with the real prover $P$ is called a simulator. The simulator does not know and cannot determine the prover's secret, and it plays the role of $P$ during the protocol.

\begin{definition}[\cite{DelfsKnebl}, Definition 4.6]
An interactive proof system $(P,V)$ is zero-knowledge if there is a probabilistic simulator $S(V^{*},x)$, running in expected polynomial time, for which every verifier $V^{*}$ outputs on input $x$ an accepting transcript $t$ of $P$ and $V^{*}$ such that these simulated transcripts are distributed in the same way as if they were generated by the honest prover $P$ and $V^{*}$.
\end{definition}
\noindent Informally, we say an interactive proof system is zero-knowledge if whatever the verifier can efficiently compute after interacting with the prover, can be efficiently simulated without interaction. The generating of transcripts includes random choices. It follows that we have a probability distribution on the set of accepting transcripts. The last condition of the definition means that the probability distribution of the transcripts which are generated by $S$ and $V^{*}$ is the same as if they were generated by the honest prover $P$ and $V^{*}$. If this distribution is not the same but statistically close to each other, then the interactive proof system is called statistical zero-knowledge. 

The statistical distance between two discrete random variables $X$ and $Y$ is defined by
\[
\Delta(X,Y) = \frac{1}{2}\sum_{k}\left|\mathbb{P}(X=k) - \mathbb{P}(Y=k)\right|, 
\]
where $k$ runs through all the objects which $X$ or $Y$ can assume and $\mathbb{P}(E)$ denote the probability of an event $E$. Let $N$ be a natural number.
We say two sets $\mathcal{X}_{N}$, $\mathcal{Y}_{N}$ of random variables are statistically close if their statitistical distance is negligible for $N \rightarrow \infty$. More precisely, if
\[
\Delta(\mathcal{X}_{N},\mathcal{Y}_{N}) = O\left(\frac{1}{p(N)}\right)
\]
for every polynomial $p(x)$ (\cite{GOL}, section 3.2.2). 

\begin{definition}[\cite{GOL}, Def. 4.3.4]
An interactive proof system $(P,V)$ is statistical zero-knowledge if there is a probabilistic simulator $S(V^{*},x)$, running in expected polynomial time, for which every verifier $V^{*}$ outputs on input $x$ an accepting transcript $t$ of $P$ and $V^{*}$ such that the distribution of these simulated transcripts is statistically close to that of they were generated by the honest prover $P$ and $V^{*}$.
\end{definition}

\section{Description of the protocol}

In this section we describe two zero knowledge protocols, which rely on the hardness of Problem \ref{IOIP}. 
\subsection{A generic protocol}

Lemma \ref{lem:np} shows that Problem \ref{IOIP} is in NP. In \cite{GMR} it is proven that every NP language admits a zero knowledge proof. In particular, this implies that one can construct an identification system whose security relies on Problem \ref{IOIP}. However, this is a purely theoretical result and the resulting protocol is inefficient. In the next subsection we propose a more direct approach which is potentially more efficient. This comes at the price that the zero knowledge property relies on a heuristic and not directly on Problem \ref{IOIP}. 


\subsection{The protocol}

In this subsection we give a high level description of the protocol and in the next subsection we provide details for key generation and interactions. 
Our protocol is based on Problem \ref{IOIP2}. 

The public key consists of two orders $\Gamma_0,\Gamma_1$ given by their regular representation (i.e., they are each given by $n^2$ matrices in $M_{n^2}(\mathbb{Q})$ which form an integral basis of the respective order). The secret key is a matrix $M\in M_{n^2}(\mathbb{Z})$ with $\Gamma_1=M^{-1}\Gamma_0M$. We denote the isomorphism corresponding to $M$ by $\phi$. Note that by a random isomorphism we mean conjugation by a suitable random integer matrix. We will specify the details of this in the next subsection. 
\begin{remark}
When $\Gamma_0$ possesses nontrivial automorphisms it might happen that $\Gamma_0=\Gamma_1$, i.e., $M$ induces an automorphism (which may be a concern since automorphisms are rare and are potentially easier to find than general isomorphisms). However, one can easily check if this is the case by deciding whether the representing $n^2\times n^2$ matrices of $\Gamma_0$ and $\Gamma_1$ generate the same $\mathbb{Z}$-lattice in $M_{n^2}(\mathbb{Q})$. In this case one should generate a new $M$. Note however, that the probability of a randomly chosen isomorphism to be an automorphism is negligible. 
\end{remark}
The steps of the protocol are as follows:
\begin{protocol}
\begin{enumerate}
    \item The prover chooses a random $r \in \{0,1\}$ and a random isomorphism $\psi:\Gamma_r\rightarrow\Gamma$ for some $\Gamma$ and sends the regular representation of $\Gamma$ corresponding to a random basis to the verifier.
    \item The verifier sends a random one-bit challenge $i\in\{0,1\}$ to the prover. 
    \item The prover computes
\[
\delta = 
\left\{
\begin{aligned}
\psi \textnormal{, if } r = i \\
\psi \circ \phi^{-1} \textnormal{, if } r = 0, i = 1 \\
\psi \circ \phi \textnormal{, if } r = 1, i = 0.
\end{aligned} \hspace*{3mm}
\right.
\]
The prover sends the isomorphism $\delta$ to the verifier.
    \item The verifier accepts if $\Gamma=\delta(\Gamma_{i})$.
\end{enumerate}
\end{protocol}
The first message is a commitment by the prover that he knows an isomorphism.
The second message is the challenge by the verifier. If the challenge sent by the verifier is the same as the prover's selection, then the prover has to open the commitment and unfold $\psi$. If not, then the verifier has to show his secret in encrypted form, by providing $\psi \circ \phi^{-1}$ or $\psi \circ \phi$.
\begin{remark}
Note that Protocol 1 is similar to Protocol 2 in \cite{GMW} based on the graph isomorphism problem. 
Furthermore, Protocol 2 in \cite{GMW} is simpler than Protocol 1 because in the first step the prover sets $r = 1$ and he does not use random selection as in Protocol 1.
\end{remark}

\subsection{Details of the protocol}

In this subsection we discuss the execution of the protocol and the key generation in more detail. 

\paragraph{Constructing orders in division algebras}

First we show how to construct division algebras of degree $p$ over $\mathbb{Q}$ where $p$ is a prime number.
\begin{enumerate}
    \item Find a cyclic Galois extension $\mathbb{L}$ of $\mathbb{Q}$ of degree $p$. Let the generator of the Galois group be $\sigma$.
    \item Find an integer $b$ which is not in the image of the norm map from
    $\mathbb{L}$ to $\mathbb{Q}$.
    \item Output the cyclic algebra $(\mathbb{L}|\mathbb{Q},\sigma,b)$.
\end{enumerate}
Now we give a brief description on how to carry out each individual step. For the first step find an integer $l$ such that $\varphi(l)$ (where $\varphi$ denotes Euler's totient fucntion) is divisible by $p$. Then the $l$th cyclotomic field $\mathbb{Q}(\epsilon_l)$ contains a subfield which is Galois over $\mathbb{Q}$ and has degree $p$ since the extension $\mathbb{Q}(\epsilon_l)|\mathbb{Q}$ is abelian and thus all its subgroups are normal subgroups. This field can be constructed by taking the fixed field of an appropriate subgroup of $Gal(\mathbb{Q}(\epsilon_l)|\mathbb{Q})$. 

The second step can be accomplished by choosing a small random $b$ (between 1 and 100 for example) and checking whether $b$ is in the image of the norm map. One has to check whether $b$ is represented locally at every prime dividing the discriminant of the norm form which is easy as $b$ is small and the discriminant of the form is small as well (this is important as one needs to factor $b$ and the discriminant of the norm form). This procedure is essentially the basis of detecting division algebras locally which is described in \cite[Section 6]{Ithesis}. The main component their is the computation of a maximal order. Once a maximal order is computed one can determine its local indices by factoring its discriminant. The probability of success is high as the image of the norm map is a subgroup of infinite index in $\mathbb{Q}$. We give an example of generating a division algebra together with an order in Section \ref{subsec:divalg}

\paragraph{Random Isomorphisms}



Now we generate a random nonsingular $n \times n$ integer matrix $M = (a_{i,j})$ in some probability distribution similarly as in \cite{Hartung}. Let $t$ be a security parameter and we choose $a_{i,j} \in (-t,t)$ uniformly at random and independently from each other for $1 < i \ne n$.
Then we compute the minors correspond to $a_{1,j}$ for every $1 \le j \le n$ (these are the entries in the first column of the adjoint matrix $M^{-1}\cdot \det(M)$). We denote these minors by $b_{1,1}, \dots{} ,b_{1,n}$. Furthermore, we can compute the entries of the first row of $M$ by solving the following linear diophantine equation:

\[
\sum_{i=1}^{n}a_{1,i}b_{1,i} = 1,
\]
which is the greatest common divisor of $b_{1,1}, \dots{} ,b_{1,n}$. We can compute the greatest common divisor by using the extended euclidean algorithm. If the greatest common divisor is not equal to $1$, then repeat the process with 
some new $a_{i,j}$, $1 < i \le n$. We compute a short solution in euclidean 2-norm of the above equation by using the LLL algorithm \cite{havas}. 
The random distribution of these matrices is denoted by $D_{t}(M)$. Formally, we have the following algorithm.


\begin{algorithm}
\caption{Computation of a random nonsingular integer matrix $M = (a_{i,j})$}
\begin{algorithmic}
\STATE choose $a_{i,j} \in (-t,t)$, $i \ne 1$ uniformly at random and independently from each other
\STATE \textbf{For} $j = 1$ to $n$
\STATE compute $b_{1,j}$, which is the minor corresponds to $a_{1,j}$.   
\STATE compute gcd($b_{1,1} \dots{} ,b_{1,n})$
\STATE \textbf{If} gcd$(b_{1,1}, \dots{} ,b_{1,n}) \ne 1$ 
\STATE repeat with a new $a_{i,j} \in (-t,t)$, $i \ne 1$. 
\STATE \textbf{Else} compute $a_{1,i}$ by solving $\sum_{i=1}^{n}a_{1,i}b_{1,i} = 1$ = gcd($b_{1,1} \dots{} ,b_{1,n})$.
\STATE reduce $(a_{1,1}, \dots{} ,a_{1,n})$ by the LLL algorithm
\RETURN $M$
\end{algorithmic}
\end{algorithm}

\begin{remark}
To find a random integer matrix we follow the algorithm of Hartung and Schnorr \cite{Hartung}. One can also
generate random isomorphisms by using matrices over $\mathbb{Q}$.
\end{remark}

The distribution of such random matrices is not necessarily uniform; actually we do not know the distribution exactly. Thus we do not know the distribution of the products and inverses of these matrices as well. We assume that these distributions statistically close. More precisely, we have the following heuristic.
\begin{heuristic}
If the random matrices $A_{1}$, $A_{2}$ $A_{3}$ and $A_{4}$ are independent and distributed according to $D_{t}(A_{i})$ $(i = 1,2,3,4$ resp.) then the distributions of $A_{1}$ and $A_{1}A_{2}$, $A_{1}A_{3}$, $A_{4}A_{1}$ and $A_{1}A_{4}^{-1}$ are statistically close to each other. 
\end{heuristic}
\begin{remark}
We have not checked this assumption experimentally yet, however, a similar assumption is made for integral quadratic forms in\cite{Hartung} and integral 
quadratic forms correspond to norm forms of orders of quaternion algebras. 
\end{remark}\

\noindent Fix an order $\Gamma_{0}$ given by its regular representation and choose a random matrix $M$ is of distribution $D_{t}(M)$ applying the above process. Compute
$\Gamma_{1} = M^{-1}\Gamma_{0}M$. The public key is $(\Gamma_{0},\Gamma_{1})$, the private key is $M$.

\paragraph{Choosing parameters} 

When choosing parameters we consider two types of attacks. One possible attack is to guess the secret isomorphism, i.e., guess the entries of the matrix $M$ as verifying whether a certain $M$ suffices can be accomplished in polynomial time. Another possible attack is to compute an isomorphism between the underlying division algebras. We stress that it is an open problem whether an isomorphism between the underlying division algebras can be applied to computing an isomorphism of orders. Nevertheless, Problem \ref{DAIP} is a potentially easier problem for which there exists an algorithm with a clear complexity analysis \cite{IRS}. 

We suggest to use orders division algebras of degree 5 and the bound for a random isomorphism to be 100. In this setting orders are represented by $25\times 25$ matrices thus and so is the secret isomorphism. This implies that searching thorough all possible $25\times 25$ matrices with entries between -100 and 100 is clearly not feasible. In order to compute an isomorphism of the underlying division algebras one has to compute an isomorphism of full matrix algebras of degree 25 \cite[Section 4]{IRS} which is infeasible as discussed in Section 2.1. Choosing the degree to be 5 comes from the fact that lower degree algebras might have better algorithms as the algorithm from \cite{IRS} (such as \cite{Pilnikova},\cite{de Graaf}).
\begin{remark}
Note that the parameters we propose allow an adversary to factor the discriminant of the orders. However, we would like to emphasize that Problem \ref{IOIP2} is still hard even when one knows the factorization of the discriminant. The most costly part of the algorithm from \cite{IRS} is the exhaustive search component not the factoring of the discriminant. 
\end{remark}

\subsection{Security of Protocol 1}
The verifier accepts the proof of a fraudulent prover with probability $\frac{1}{2}$. Iterating the protocol
$k$ times independently and sequentially, the probability of the cheating can be reduced to $2^{-k}$.
On the other hand, the only information a honest prover provides to the verifier is the fact that he knows an isomorphism between the two orders.

\noindent \textbf{Completeness}
It is clear that if $r = i$, then $\delta(\Gamma_{i}) = \psi(\Gamma_{i}) = \Gamma$ and if $N\in M_{n^2}(\mathbb{Z})$ is the matrix corresponds to $\psi$, then $N^{-1}\Gamma_iN = \Gamma$ and the verifier will accept.
If $r = 0$, $i = 1$, then $(\psi \circ \phi^{-1})(\Gamma_{1}) = \psi(\phi^{-1}(\Gamma_{1})) = \psi(\Gamma_{0}) = \Gamma$ for some $\Gamma$. More precisely, $(M^{-1}N)^{-1}\Gamma_{1}(M^{-1}N) = N^{-1}(M\Gamma_{1}M^{-1})N = N^{-1}\Gamma_{0}N = \Gamma$, thus the verifier will accept. Finally, if  $r = 1$, $i = 0$, then $(\psi \circ \phi)(\Gamma_{0}) = \psi(\phi(\Gamma_{0})) = \psi(\Gamma_{1}) = \Gamma$ for some $\Gamma$. More precisely, $(MN)^{-1}\Gamma_{0}(MN) = N^{-1}(M^{-1}\Gamma_{0}M)N = N^{-1}\Gamma_{1}N = \Gamma$, thus the verifier will accept. 
These facts impliy that if the prover knows $\phi$, and both the prover and the verifier follow the protocol, then 
the verifier will always accept.


\noindent \textbf{Soundness}
We prove that if $P^{*}$ is a fraudulent prover, then the verifier $V$ will reject with probability at least $\frac{1}{2}$. If any $P^{*}$ can convince the verifier with both challenges $i = 0, 1$, then clearly there exist $N_{0}, N_{1}\in M_{n^2}(\mathbb{Z})$
such that $N_{0}^{-1}\Gamma_{0}N_{0} = \Gamma$ and $N_{1}^{-1}\Gamma_{1}N_{1} = \Gamma$.
Then we have $N_{0}^{-1}\Gamma_{0}N_{0} = N_{1}^{-1}\Gamma_{1}N_{1}$ and so
$\Gamma_{1} = (N_{1}^{-1}N_{0})^{-1}\Gamma_{0}(N_{1}^{-1}N_{0})$ which gives another integral isomorphism between $\Gamma_{0}$ and $\Gamma_{1}$ i.e., another private key.
It follows that at most one of the challenges may lead to acceptance. Hence, with probability at least $\frac{1}{2}$, the verifier will then reject.

\begin{proposition}
Under Heuristic 1, Protocol 1 is statistical zero knowledge.
\end{proposition}

\begin{proof}
According to Theorem 2 in \cite{GAV} and in \cite{HRV}, it is enough to prove the proposition for honest verifier.
Let $O$ 
denote the set of orders given by their regular representations and let $H$ denote the set of isomorphisms from $\Gamma_{0}$ to $\Gamma$ and from $\Gamma_{1}$ to $\Gamma$.
The set of accepting transcripts is
\[
\{(\Gamma, i, \delta) \in O \times \{0,1\} \times H: \delta(\Gamma_{i}) = \Gamma\}.
\]
We describe a simulator $S$, which satisfies the desired properties.
\begin{algorithm}
\caption{Simulator $S$}
\begin{algorithmic}
\STATE transcript $S$(algorithm of $V$, $\Gamma_{0}$ and $\Gamma_{1}$ given by regular representations)
\STATE choose a random $s \in \{0,1\}$ and an isomorphism $\hat{\delta}$ from $\Gamma_{s}$ at random according to $D_{t}(\Gamma_{s})$. 
\STATE $\hat{\Gamma} \leftarrow \hat{\delta}(\Gamma_{s})$
\STATE choose an isomorphism $\eta$ from $\hat{\Gamma}$ at random according to $D_{t}(\hat{\Gamma})$.
\STATE $\Gamma^{'} \leftarrow \eta(\hat{\Gamma})$
\STATE choose an isomorphism $\theta$ from $\hat{\Gamma}$ at random according to $D_{t}(\hat{\Gamma})$.
\STATE $\Gamma^{''} \leftarrow \theta(\hat{\Gamma})$
\STATE$i \leftarrow V(\Gamma^{'})$
\STATE \textbf{If} $s = i$
\RETURN $(\Gamma^{'}, s, \eta \circ \hat{\delta})$
\STATE \textbf{If} $s \ne i$
\RETURN $(\Gamma^{''}, s, \theta \circ \hat{\delta})$
\end{algorithmic}
\end{algorithm}


\noindent The simulator $S$ uses the verifier $V$ to get the challenge $i$ and it tries to find out $i$ in advance. If $S$ was successful in guessing $i$, he can provide a valid transcript $(\hat{\Gamma}, s, \eta \circ \hat{\delta})$. 
Now we prove that under Heuristic 1, the distribution of the simulator's output is statistically close to an output coming 
from an interaction between a honest prover $P$ and an honest verifier $V$. It is clear from the definition of $S$ that
the distribution of $s$ is the same as the challenge $i$ coming from the interaction between $P$ and $V$ in the protocol.
Define $K$, $N_{1}$, $N_{2}$ by $K^{-1}\Gamma_{s}K = \hat{\Gamma}$,  $N_{1}^{-1}\hat{\Gamma} N_{1} = \Gamma^{'}$, $N_{2}^{-1}\hat{\Gamma} N_{2} = \Gamma^{''}$, respectively. If $i = s$, then the output of $P$ is $K \in D_{t}(K)$ and the output of $S$ is $KN_{1}$, where $N_{1} \in D_{t}(N_{1})$.  
According to Heuristic 1, the distribution of $K$ and $KN_{1}$ are statistically close to each other. 
If $i \ne s$, then the output of $P$ is $KM^{-1}$ or $MK$ and the output of $S$ is $KN_{2}$, where
$M \in D_{t}(M)$, $K \in D_{t}(K)$ and $N_{2} \in D_{t}(N_{2})$.
According to Heuristic 1, the distribution of $KM^{-1}$, $MK$ and $KN_{2}$ are statistically close to each other. Thus the outputs of $P$ and $S$ are statistically close and so Protocol 1 is statistically zero knowledge under Heuristic 1.


\end{proof}
\begin{remark}\label{remark}
The proof of the zero knowledge property of the protocol is reminiscent of the zero knowledge property of the classical interactive proof system based on graph isomorphisms. However, in the case of graph isomorphisms, if you consider the product of two random isomorphisms it is indistinguishable from a random isomorphism chosen from the same distribution (i.e., selecting a unformly random element of the symmetric group). However, if one looks at the composition of two order isomorphisms, the corresponding conjugating matrices multiply which will have a different distribution. This is why we have to modify the simulator and can only claim statistical zero knowledge similarly to \cite{Hartung}. 
\end{remark}

Finally, we briefly comment on the impact of the proposed protocol. Problem \ref{IOIP2} is related to well-studied number theoretical problems and is somewhat connected to computational assumptions reminiscent of lattice-based and multivariate assumptions. This provides some motivation that Problem \ref{IOIP2} is hard even for a quantum computer. Furthermore, if one could extend our protocol to be able to handle arbitrarily large challenges, one could build digital signature schemes which could have competitive signing and verification speed. The reason for this is that by choosing a larger division algebra, one could potentially choose the conjugating matrix $M$ to have small entries which enables fast matrix multiplication. Note that the key sizes will still be very large. Thus we believe that this paper could potentially be a starting point for a new line of post-quantum schemes. 
\section{A toy example}\label{sec:example}
\subsection{Identification scheme}
In this subsection we give a concrete example of our identification scheme. This is just a small example meant to provide some clarity, so the parameters used are not meant to provide sufficient security. 

Let $\A$ be the quaternion algebra with quaternion basis $1,u,v,uv$ where $u^2=-1,~v^2=3,~uv+vu=0$. Let $O_1$ be the $\mathbb{Z}$-lattice generated by $1,u,v,uv$. It is clear that $O_1$ is also a ring, thus it is an order in $A$. First we compute the regular representation of $O_1$:
$$
\begin{pmatrix}
1&0&0&0 \\
0&1&0&0 \\
0&0&1&0 \\
0&0&0&1
\end{pmatrix},
\begin{pmatrix}
0 &-1 & 0 & 0\\
1& 0 &0 & 0  \\
0&  0&  0& -1\\
0&  0&  1&  0
\end{pmatrix}, 
\begin{pmatrix}
 0 & 0 & 3 & 0 \\
 0 & 0 & 0 &-3 \\
 1 & 0 & 0 & 0 \\
 0 &-1 &  0 &  0 
\end{pmatrix},
\begin{pmatrix}
0& 0& 0& 3 \\
0& 0& 3& 0 \\
0& 1& 0& 0 \\
1& 0& 0& 0
\end{pmatrix}.
$$
Now we construct an order $O_2$ isomorphic to $O_1$ by first choosing a random matrix $B$ and computing a random basis of $O_2=B^{-1}O_1B$. Here we choose the following matrix:
$$
B=\begin{pmatrix}
-1 &-5& -1& -5 \\
10 & 3& 14& 14 \\
10&  4& 15&  0 \\
 7&  9&  9& 18
\end{pmatrix}.
$$
The regular representation of $O_2$ is given by the following matrices:
$$
\begin{pmatrix}
 22823&  25756&  24424& 127163 \\
  -426&   -488&   -453&  -2410 \\
-15106& -17043& -16168& -84136 \\
 -1110&  -1252&  -1188&  -6183 
\end{pmatrix}, 
\begin{pmatrix}
  -236067 & -749009&  -241581& -1276977 \\
    4408  &  14188 &    4479 &   24201 \\
  156222 &  495551 &  159890 &  844853 \\
   11510&    36424 &   11794 &   62093 
\end{pmatrix},
$$
$$
\begin{pmatrix}
 -700782 &-1036680  &-856635 &-2353415 \\
   13175  &  19601   & 16073  &  44639 \\
  463691  & 685895   &566827  &1557060 \\
   34088  &  50397    &41679  & 114354  
\end{pmatrix},
\begin{pmatrix}
 -434287  &-941658 &  -509832  &-1628214 \\
    8154  &  17819  &   9546  &  30876 \\
  287358 &  623010  & 337355 & 1077234 \\
   21144 &    45786  &  24834 &   79133
\end{pmatrix}.
$$
The regular representation of $O_1$ and $O_2$ is part of the public key, the matrix $B$ is secret. 
Now the prover chooses a random matrix $C$ and computes $O_3=C^{-1}O_2C$. Here we choose 
$$
C=\begin{pmatrix}
 -4& -1& -1& -4 \\
13 & 9 & 3  &6 \\
 2 & 9 &15 &13 \\
 3 &13 &12  &5
\end{pmatrix}.
$$ 
The the prover sends over a random basis of $O_3$, in our case this will be the following: 
$$
\begin{pmatrix}
 308953146301 & 655710307003 & 579155375644 & 361329662174 \\
-383835331750 &-814637386699 &-719527537820 &-448906551540 \\
 476972774120 &1012308722844 & 894120520149 & 557833493285 \\
-332130522360 &-704901082168 &-622603073690 &-388436279011
\end{pmatrix},
$$
$$
\begin{pmatrix}
  587711341384 & 1356738677803 & 1237179068072  & 755461309928 \\
 -730157244630 &-1685576756090 &-1537039017600  &-938565435976 \\
907329516990  &2094581071851 & 1910000731132  &1166307849683 \\
 -631800895590 &-1458519943242 &-1329991087710  &-812135315774
\end{pmatrix},
$$

$$
\begin{pmatrix}
  735180739907 & 1574065359241 & 1394017791570 &  866716597474 \\
 -913369378520 &-1955577756259 &-1731891356960 &-1076786103514 \\
 1134997980115 & 2430097676536 & 2152133892443 & 1338067692447 \\
 -790333311240 &-1692150272837 &-1498595710360  &-931736874191 
\end{pmatrix},
$$
$$
\begin{pmatrix}
 981143888960 & 2111674730047 & 1869491446440 & 1157041030470 \\
-1218947580410 &-2623489619892 &-2322607428900 &-1437477610080 \\
 1514724572734 & 3260078004908 & 2886186911608 & 1786280805793 \\
-1054748385480 &-2270090596754 &-2009738956870 &-1243841177300
\end{pmatrix}.
$$
\begin{remark}
Note that the matrices corresponding to $O_3$ have larger coefficients than the matrices corresponding to $O_2$. This is an example of the phenomenon explained 
at the end of Remark \ref{remark}.
\end{remark}
Now the prover chooses a random bit $b$. If $b=0$, then the prover reveals the matrix $C$. In this case the verifier first computes $C^{-1}O_2C$ (i.e., conjugates the given four matrices by $C$). Then the verifier accepts if the given four matrices of generate the same $\mathbb{Z}$-lattice as the given four matrices (i.e., checks that the transition matrix is an integer matrix with determinant $1$ or $-1$). If $b=1$, then the prover reveals $BC$ and the verifier checks in a similar fashion. 
\subsection{Generating division algebras and orders}\label{subsec:divalg}

We provide an example of how to generate division algebras and orders in them. The outline will be the following. We find a cyclic extension $\mathbb{K}$ of $\mathbb{Q}$ and then an element $b\in\mathbb{Q}$ for which the cyclic algebra $\A=(\mathbb{K}|\mathbb{Q},\sigma,b)$ is a division algebra. Then we provide an order in $\A$. 

In order to make our life easier we will look for a cyclic algebra of degree 5, i.e., a field $\mathbb{K}$ which is a Galois extension of $\mathbb{Q}$ is of degree 5. Since 5 is a prime number, $\A=(\mathbb{K}|\mathbb{Q},\sigma,b)$ is a division algebra if and only if $b$ is not a norm in the extension $\mathbb{K}|\mathbb{Q}$. We are looking for $\mathbb{K}$ as a subfield of a certain cyclotomic field because cyclotomic extensions are Abelian hence every subfield of them is automatically a Galois extension of $\mathbb{Q}$. 

Since $\varphi(11)=10$, the 11th cyclotomic field contains a subfield $\mathbb{K}$ of degree 5 over $\mathbb{Q}$ which is the splitting field of the polynomial $x^5-11x^4+44x^3-77x^2+55x-11$. The Galois group of $\mathbb{K}$ over $\mathbb{Q}$ is cyclic because every group of order 5 is cyclic. Now we need to find a $b$ which is not a norm in the extension $\mathbb{K}|\mathbb{Q}$. Checking that a certain $b$ is a norm or not can be accomplished efficiently using \cite[Section 6]{Ithesis} if one can factor $b$. In this setting $b$ does not need to be large and most $b$-s aren't norms so this can actually be done by guessing easily. Using the computational algebra system MAGMA \cite{Magma} one can compute that $b=2$ is not a norm. The discussion so far implies that $\A=(\mathbb{K}|\mathbb{Q},\sigma, 2)$ is a division algebra. 

Let $u\in \A$ be such that $u^5=2$ (which is given by the definition of a cyclic algebra). Then the set $\sum_{i=0}^4 Ou^i$ where $O$ is the ring of integers of $\mathbb{K}$ is an order. Indeed, it is full $\mathbb{Z}$-lattice and a subring of $\A$ which contains 1. Alternatively, one can also generate an order by selecting a $\mathbb{Q}$-basis of $\A$ containing 1 and multiply every basis element with a suitable integer to make structure constants integral. 
\begin{remark}
This method easily generalizes to constructing division algebras of square-free degree. Indeed, the tensor product of division algebras of coprime degrees is again a division algebra \cite{Pi}.
\end{remark}
\begin{remark}
We do not see any security risk in setting this exact division algebra as a global parameter (i.e., this division algebra can be used in any protocol execution). 
\end{remark}

\paragraph*{Acknowledgement.}
We would like to thank the anonymous reviewers for the careful reading and the helpful suggestions which have improved the quality of this paper considerably. 
S\'andor Z. Kiss was supported by the Hungarian National Research, Development and Innovation Office - NKFIH, Grants No. K109789, K129335, K115288. Sándor Z. Kiss was supported by the János Bolyai Research Scholarship of the Hungarian 
Academy of Sciences and by the \'UNKP-18-4 New National Excellence Program of the Ministry of Human Capacities. Sándor Z. Kiss was supported by the \'UNKP-19-4 New National Excellence Program of the Ministry for Innovation and Technology. Supported by the \'UNKP-20-5 New National Excellence Program of the Ministry for Innovation and Technology from the source of the National Research, Development and Innovation Fund. P\'eter Kutas was supported by an EPSRC New Investigator grant (EP/S01361X/1).

\end{document}